\documentclass[amstex,12pt]{article}

\usepackage{graphicx}
\usepackage{amsfonts}
\usepackage{amssymb}
\usepackage{amsmath}
\usepackage{latexsym}
\usepackage{enumerate}

\newtheorem{thm}{Theorem}[section]
\newtheorem{lm}[thm]{Lemma}
\newtheorem{prop}[thm]{Proposition}

\newtheorem{q}[thm]{Question}

\newtheorem{ex}[thm]{Example}
\newtheorem{df}[thm]{Definition}

\begin{document}

\author{Jin-ichi Itoh and Costin V\^\i lcu}
\title{Every graph is a cut locus}
\maketitle

\noindent {\bf Abstract} {\small
We show that every connected graph can be realized as the cut locus of some point on some Riemannian surface $S$ which, in some cases, has constant curvature. We study the stability of such realizations, and their generic behavior.
\\{\bf Math. Subj. Classification (2010):} 53C22, 05C62}


\section{Introduction}

Unless explicitly stated otherwise, by a Riemannian manifold here we always mean a complete, 
compact and connected manifold without boundary. 
We shall work most of the time with surfaces ($2$-dimensional manifolds) $S$, 
and let $M$ denote manifolds of arbitrary dimension $d$.

All graphs we consider in the following are finite, connected and may have loops and multiple edges. 
For the simplicity of our exposition, we see every graph $G$ as a $1$-dimensional simplicial complex.
The {\it cyclic part} of $G$ is the minimal (with respect to inclusion) subset $G^{cp}$ of $G$, 
to which $G$ is contractible; i.e., $G^{cp}$ is the minimal subset of $G$ obtained by repeatedly contracting edges incident to degree one vertices, and for each remaining vertex of degree two (if any) merging its incident edges. $G^{cp}$ thus inherits a natural structure of simplicial complex.
A graph is called {\it cyclic} if it is equal to its cyclic part, 
and it is called {\it regular} if all its vertices have the same degree.
A {\it length graph} is a weighted graph with positive weights; i.e., each edge is endowed with a positive number (also called length).

The notion of cut locus was introduced by H. Poincar\'e \cite{p} in 1905, 
and gained since then an important place in global Riemannian geometry.
The {\it cut locus $C(x)$ of the point} $x$ in the Riemannian manifold $M$ is the set of all extremities 
(different from $x$) of maximal (with respect to inclusion) {\it segments} (i.e., shortest geodesics) starting at $x$; 
for basic properties and equivalent definitions refer, for example, to \cite{ko} or \cite{sa}.

For Riemannian surfaces $S$ is known that $C(x)$, if not a single point, is a local tree 
(i.e., each of its points $z$ has a neighborhood $V$ in $S$ such that 
the component $K_z(V)$ of $z$ in $C(x)\cap V$ is a tree), 
even a tree if $S$ is homeomorphic to the sphere.
A {\it tree} is a set $T$ any two points of which can be joined
by a unique Jordan arc included in $T$.
The {\it degree} of a point $y$ of a local tree is the number of components
of $K_y(V)\setminus \{y\}$ if $V$ is chosen such that $K_y(V)$ is a tree.

S. B. Myers \cite{M} for $d=2$, and M. Buchner \cite{Bu2} for general $d$, established that 
the cut locus of a real analytic Riemannian manifold of dimension $d$ is homeomorphic to 
a finite $(d-1)$-dimensional simplicial complex.
For a class of Liouville manifolds, in particular for hyperellipsoids in the Euclidean space $\mathbb{R} ^d$, 
the cut locus is reduced to a disc of dimension of most $(d-1)$, see \cite{Itoh-Kiyohara_1} and \cite{Itoh-Kiyohara_2}.

For Riemannian metrics on $S$ non-analytic, cut loci may be quite large sets.
J. Hebda \cite{H1} showed, for any ${\cal C}^\infty$ metric on $S$, that 
the Hausdorff 1-measure of any compact subset of the cut locus of any point is finite. 
Independently and using different techniques, 
J. Itoh \cite{I2} proved the same result under the weaker assumption of ${\cal C}^2$ metric. 
The differentiability of the metric cannot be lowered more;
for example, the main result in \cite{z-ep} states that 
on most (in the sense of Baire category) convex surfaces (known to be of  differentiability class ${\cal C}^1 \setminus {\cal C}^2$), 
most points are endpoints of any cut locus.

The problem of constructing a Riemannian metric with preassigned cut locus on a given manifold 
also received a certain interest. 
H. Gluck and D. Singer \cite{GS} constructed a Riemannian metric such that 
a non triangulable set, consisting of infinitely many arcs with a common extremity, becomes a cut locus. 
Another example of infinite length cut locus was provided by J. Hebda \cite{H2}, 
while the case of a submanifold as preassigned cut locus was considered by L. B\'erard-Bergery \cite{Bb}.
J. Itoh \cite{I1} showed that for any Morse function on a differentiable surface $S$, 
with only one critical point of index $0$ and no saddle connection, 
there exists a Riemannian metric on $S$ with respect to which $C_f$, 
the union of all unstable manifolds of critical points of $f$ with positive index, becomes a cut locus.
Independently but in the same direction as \cite{I1}, M. Y. Park showed that,
under some sufficient conditions, for any smoothly embedded, connected, finite cubic graph $G$ in the surface $S$, there exists a Riemannian metric $\alpha$ on $S$ and a point $x$ in $S$ such that the cut locus of $x$ with respect to $\alpha$ is $G$ \cite{Park1}, and that this cut locus is stable with respect to the metric \cite{Park2}.
All these results assume the manifold be given, 
and search for a metric with respect to which some subset of the manifold becomes a cut locus.

A different approach was considered in \cite{iv1}, where the authors showed that any combinatorial type of finite tree can be realized as a cut locus on some, initially unknown, doubly covered convex polygon.

Our results here give that approach much more generality, by showing (see Theorem \ref{real_cl}) that every connected length graph can be realized as a cut locus; 
i.e., there exist a Riemannian surface $S_G=(S_G,h)$ and a point $x \in S_G$ such that $C(x)$ is isomorphic to $G$.
This is a partial converse to Myers' theorem mentioned above.
If moreover $G$ is cyclic and regular then it can be realized on a surface of constant curvature (Theorem \ref{CL_ct_curv}).
In the second part of this paper we show that --roughly speaking-- stability is a generic property of cut locus realizations.

In a forthcoming paper \cite{iv3} we are concerned about 
the orientability of the surfaces $S_G$ realizing the graph $G$ as a cut locus.

Employing the notion of {\sl cut locus structure} \cite{iv2}, one can also regard our results as completing with additional information the surface case in the results of Buchner \cite{Bu1}, \cite{Bu2}, \cite{Bu3}.

Recently, and from a viewpoint different from ours, cut loci and infinite graphs were studied by O. Baues and N. Peyerimhoff \cite{Ba-1}, \cite{Ba-2}, and by M. Keller \cite{Ke}, while in discrete group theory a similar notion, {\sl dead-end depth}, was studied by S. Cleary and T. R. Riley \cite{CR}, and by T. R. Riley and A. D. Warshall \cite{RW}.


\section{Every graph is a cut locus}
\label{every_graph}

Recall that {\sl a segment between a point $x$ and a closed set} $K$ not containing $x$ is 
a segment from $x$ to a point in $K$, not longer than any other such segment; 
the {\sl cut locus $C(K)$ of the closed set} $K \subset S$ is 
the set of all points $y \in S$ such that there is a segment from $y$ to $K$ not extendable as a segment beyond $y$.

\begin{df}
{\rm Let $G$ be a graph. A }{\it strip} {\rm on $G$ (in short, a $G$-}{\sl strip}{\rm ) is a topological surface $P_G$ with boundary, such that:}
\\{\sl i)} {\rm the boundary of $P_G$ is homeomorphic to a circle, and}
\\{\sl ii)} {\rm  $P_G$ contains (a graph isomorphic to) $G$ and is contractible to $G$.

A }{\it Riemannian} $G$-{\it strip} {\rm is a $G$-strip $P_G$ endowed with a Riemannian metric such that 
the cut locus of ${\rm bd}(P_G)$ in $P_G$ is precisely $G$.
If $G$ is a length graph, we ask in addition that the induced lengths on the edges of $G$ by the metric of $P_G$ 
coincide to the original weights.}
\end{df}

Regarding condition {\sl (ii)} above, the simple example of a tree on a cylinder shows that a topological surface with boundary is not contractible to each graph it contains.

\begin{df}
{\rm We say that a graph (or a length graph)} $G$ {\it can be realized as a cut locus} {\rm if 
there exist a Riemannian surface $S_G=(S_G,h)$ and a point $x$ in $S_G$ such that $G$ is isomorphic to $C(x)$.}
\end{df}

A. D. Weinstein (Proposition C in \cite{W}) proved the following.

\begin{lm}
\label{Wein}
Let $M$ be a $d$-dimensional Riemannian manifold and $D$ an $d$-disc embedded in $M$.
There exists a new metric on $M$ agreeing with the original metric on a neighborhood 
of $M \setminus ({\rm interior \; of \;} D)$ such that, for some point $p$ in $D$, 
the exponential mapping at $p$ is a diffeomorphism 
of the unit disc about the origin in the tangent space at $p$ to $M$, onto $D$.
\end{lm}

\begin{prop}
\label{glue}
The following statements are equivalent:

i) the length graph $G$ can be realized as a cut locus;

ii) there exists a $G$-strip;

iii) there exists a Riemannian $G$-strip.
\end{prop}

\noindent{\sl Proof:} {\it (i) $\to$ (ii)} 
Consider a point $x$ on a Riemannian surface $(S,g)$, and a segment $\gamma: [0,l_\gamma] \to S$ parametrized by arclength, with $\gamma(0)=x$ and $\gamma(l_\gamma) \in C(x)$.
For $\varepsilon >0$ strictly smaller than the injectivity radius ${\rm inj}(x)$ at $x$,
the point $\gamma(l_\gamma-\varepsilon)$ is well defined because ${\rm inj}(x) \leq l_\gamma$.
Since $S \setminus C(x)$ is contractible to $x$ along geodesic segments, and thus homeomorphic to an open disk, the union over all segments $\gamma$ of those points $\gamma(l_\gamma-\varepsilon)$ is homeomorphic to the unit circle.

{\it (ii) $\to$ (iii)} An explicit construction of a Riemannian $G$-strip from a given $G$-strip was provided by 
the first author in \cite{I1}.

{\it (iii) $\to$ (i)} A. D. Weinstein's result above (Lemma \ref{Wein}) shows that, given a Riemannian $G$-strip $P_G$, 
one can glue it to a disk to obtain a surface $S_G$, and there exists a metric $g$ on $S_G$ 
agreeing with the original metric on $P_G$, and a point $x$ in $S_G$ with $C(x)=G$.
\hfill $\Box$

\bigskip

We need one more result, well known in the graph theory.

\begin{lm}
\label{qmn}
For every graph with $m$ edges, $n$ vertices, and $q$ generating cycles holds $q=m-n+1$.
\end{lm}

\begin{thm}
\label{real_cl}
Every length graph can be realized as a cut locus.
\end{thm}

\noindent{\sl Proof:}
By Proposition \ref{glue}, it suffices to provide, for every length graph $G$, at least one $G$-strip.

We notice first that we can reduce our problem to the cyclic part $G^{cp}$ of $G$.
Assume $G \setminus G^{cp}$ consists of finitely many finite trees, say $T_1,T_2,...,T_m$. 
Since every tree $T$ has a ``leaf''-type $T$-strip, one can attach (in a natural way) 
all the $T_i$-strips to a $G^{cp}$-strip to obtain a $G$-strip.

We proceed by induction over the number $k$ of generating cycles of $G$. 

For $k=0$ and $G=G^{cp}$ the strip is elementary.

For $k=1$ and $G=G^{cp}$ our strip is the flat compact M\"obius band.

Assume now that there exist strips for all graphs with $k$ generating cycles, for some $k \geq 1$.

Let $G_{k+1}=G_{k+1}^{cp}$ be a length graph with $k+1$ generating cycles, 
and $e$ an edge of $G_{k+1}$ in some generating cycle of $G_{k+1}$.

\begin{figure*}[ht]
\centering
  \includegraphics[width=1.0\textwidth]{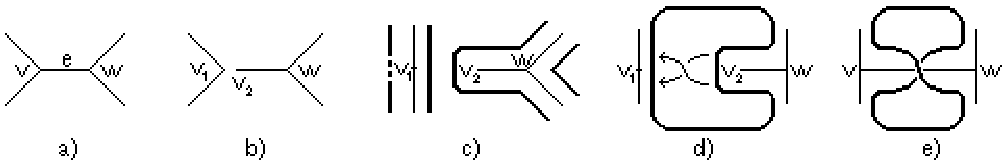}
\caption{Induction reduction: edge $e$ joins distinct vertices $v \neq w$.}
\label{Fig_1}
\end{figure*}

\begin{figure*}[ht]
\centering
  \includegraphics[width=1.0\textwidth]{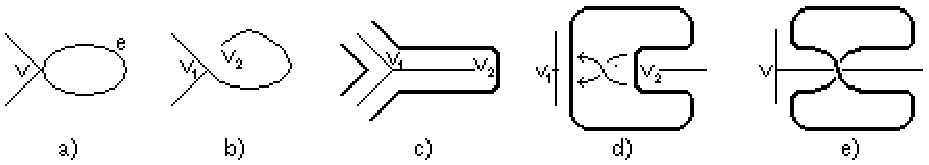}
\caption{Induction reduction: edge $e$ is a loop at $v$.}
\label{Fig_2}
\end{figure*}

Detach $e$ from $G_{k+1}$ at one extremity, say $v$; 
Figure \ref{Fig_1}(a)-(b) presents the case when $e$ joins distinct vertices $v \neq w$, 
while Figure \ref{Fig_2}(a)-(b) presents the case $v=w$.
Denote by $G_k$ the resulting length graph, and by $v_1, v_2$ the images of $v$ in $G_k$.

Since $G_k$ has one vertex more than $G_{k+1}$, it has $k$ generating cycles (see Lemma \ref{qmn}), 
and by the induction assumption there exists a $G_k$-strip $P_{G_k}$ (see Figures \ref{Fig_1}(c) and \ref{Fig_2}(c)).
Consider a planar representation of the boundary of $P_{G_k}$ as a simple closed curve
(illustrated in Figures \ref{Fig_1}(d) and \ref{Fig_2}(d)), 
and attach to it a switched $e$-strip, see Figures \ref{Fig_1}(e) and \ref{Fig_2}(e), to obtain a $G_{k+1}$-strip.
\hfill $\Box$

\bigskip

Disconnected graphs can as well be realized as cut loci, but on non-complete surfaces. To see this, consider a disconnected graph $G'$ as a subgraph of a connected graph $G$, and realize $G$ as a cut locus on a surface $S$; i.e., $G=C(x)$ for some point $x$ in $S$. With $T=G \setminus G'$, we have $C(x)=G'$ on $S \setminus T$.

\bigskip

Our Theorem \ref{real_cl} shows, in particular, that for every connected graph $G$,
there exists a $2$-cell embedding with just one face, onto some surface $S_G$.
This result is well known in the topological graph theory, see e.g. \cite{Mo-To}.

\begin{q} 
{\rm Several open questions naturally arise from Theorem \ref{real_cl}.}

i) {\rm Can the metric of the surfaces $S_G$, realizing $G$ as a cut locus, be chosen analytic?
See the result of S. B. Myers {\rm \cite{M}} mentioned in the introduction.}

ii) {\rm Cut loci on Riemannian surfaces may be quite large sets, see the introduction.
Can Theorem \ref{real_cl} be extended to infinite graphs?}

iii) {\rm Can Theorem \ref{real_cl} be extended to higher dimensions?}
\end{q}


There usually are many strips on the same graph; we formalized this by several concepts \cite{iv1}, that we briefly present next.

\begin{df}
{\rm A} {\it cut locus structure} {\rm (in short, a} {\it CL-structure}) {\rm on the graph $G$ is a strip on the cyclic part $G^{cp}$ of $G$.}
\end{df}

\begin{df}
\label{clns}
{\rm Consider, for a point $x$ on a Riemannian surface $(S,g)$ and for $\varepsilon >0$ small enough, 
the $C(x)$-strip obtained as the union, over all segments $\gamma$ starting at $x$ and parametrized by arclength, of the points $\gamma(l_\gamma - \varepsilon)$.
We call the CL-structure constructed in this way the} {\it cut locus natural structure} {\rm defined by $x$, 
and denote it by $CLNS(x)$, or by $CLNS(x, g)$ if to point out (the dependence on) the metric $g$.}
\end{df}

With these notions, Theorem \ref{real_cl} can be rephrased as that {\sl each graph possess at least one CL-structure}, while Proposition \ref{glue} and Lemma \ref{Wein} say that {\sl each CL-structure can be realized in a natural way}.

\bigskip

In order to easier handle a CL-structure, we associate to it an object of combinatorial nature.

An {\sl elementary strip} of a G-strip $P_G$ is an {\sl edge-strip} (a strip on an edge of $G$) or a {\sl point-strip} (a strip on a vertex of $G$), included in  $P_G$.
So we can think about a G-strip as union of elementary strips corresponding to all edges and vertices in $G$.
Denote by ${\cal P}$ and ${\cal A}$ the set of the point-strips, respectively edge-strips, of a CL-structure ${\cal C}$ on the graph $G$.

Below, $V$ denotes the vertex set of $G$, $E$ the edge set of $G$, while
$\bar 0$ and $\bar 1$ are the elements of the $2$-element group $(\mathbb{Z}_2, \oplus)$.

\begin{df}
{\rm Consider a $G$-strip $P_G$ as union of elementary strips, each of which has a distinguished face labeled $\bar 0$.
The face opposite to the distinguished face will be labeled $\bar 1$. 

To each pair $(v,e) \in V \times E$ consisting of a vertex $v$ and an edge $e$ incident to $v$, 
we associate the $\mathbb{Z}_2$-sum $\bar s (v,e)$ of the labels of the elementary strips $\nu \in {\cal P}$, 
$\varepsilon \in {\cal A}$ associated to $v$ and $e$; i.e.,
$\bar s (v,e) = \bar 0$ if the distinguished faces of $\nu$ and $\varepsilon$ agree to each other, and $\bar 1$ otherwise.
Therefore, to any cut locus structure ${\cal C}$ we can associate a function $s_{\cal C} :E \to \mathbb{Z}_2$, 
\begin{eqnarray}
\label{companion}
s_{\cal C}(e)= \bar s (v,e) \oplus \bar s (v',e),
\end{eqnarray}
where $v$ and $v'$ are the vertices of the edge $e \in E$.

We call the function $s_{\cal C}$ defined by (\ref{companion}) the} {\it companion function} {\rm of ${\cal C}$. }
\end{df}

\begin{df}
\label{diff_CL}
{\rm Consider two CL-structures ${\cal C}$, ${\cal C}'$ on the graph $G$.

The companion functions $s_{\cal C}$ and $s_{{\cal C}'}$ are called} {\it equivalent} {\rm on a $2$-connected component $G_{2c}$ of $G$ if they are equal, up to a simultaneous change of the distinguished face for all elementary strips in $G_{2c}$: either $s_{\cal C} = s_{{\cal C}'}$, or $s_{\cal C} =\bar 1 \oplus s_{{\cal C}'}$, on $G_{2c}$.

${\cal C}$ and ${\cal C}'$  are called} {\it equivalent} {\rm if their companion functions are equivalent on every $2$-connected component of $G$.}
\end{df}

The next sections are related to the following.

\begin{q}
\label{all_x_CL}
{\rm What can be said about the Riemannian surface $S$ if $CLNS(x)$ and $CLNS(y)$ are equivalent, for any points $x,y \in S$?}
\end{q}

From now on, all CL-structures will be considered up to equivalence.
This will allow us, whenever we consider surfaces realizing the graph $G$ as a cut locus, 
to actually think about CL-structures and their companion functions on $G$.


\section{Constant curvature realizations}
\label{ct_curv}

In this short section we present a direct way to realize some graphs as cut loci, different from that provided by Theorem \ref{real_cl}.

\begin{thm}
\label{CL_ct_curv}
Every CL-structure on a regular graph can be realized on a surface of constant curvature.
\end{thm}

\noindent{\sl Proof:}
Denote by $G$ a $k$-regular cyclic graph, and by ${\cal C}$ a CL-structure on $G$.

If $G$ is a point then the unique CL-structure on $G$ can be realized as $CLNS(x)$ 
for any point $x$ on the unit $2$-dimensional sphere.

Assume now that $G$ is a cycle. Then again we have a unique CL-structure on $G$, 
and it can be realized as $CLNS(x)$ for any point $x$ on the standard projective plane.

Consider now a graph $G$ with $q \geq 2$ generating cycles; by Lemma \ref{qmn}, we get $m\geq 2$.

For $m=2$, let $F_{2m}=F_4$ denote the square in the Euclidean plane $\Pi$.

For $m=3$, let $F_{2m}=F_6$ denote the regular hexagon in $\Pi$.

For $m \geq 4$, consider a regular $2m$-gon $F_{2m}=\bar z_1...\bar z_{2m}$ 
in the hyperbolic plane $\mathbb{H}^2$ of constant curvature $-1$, 
such that $\angle \bar z_i \bar z_{i+1} \bar z_{i+2} =2\pi/k$ (all indices are taken (mod $2m$)).

We view now the CL-structure ${\cal C}$ on $G$ as a closed path $D$ in $G$ containing all edges of $G$ precisely twice,
hence every vertex of $G$ appears precisely $k$ times in $D$.

We identify now the path $D$ with (the boundary of) $F_{2m}$, such that 
each image in $D$ of an edge of $G$ corresponds to precisely one edge in $F_{2m}$, 
each image in $D$ of a vertex of $G$ corresponds to precisely one vertex in $F_{2m}$, 
and the order of edges and vertices along $D$ is preserved.
It remains to identify, for every edge $e$ in $G$, its two images in $F_{2m}$, 
to obtain a differentiable surface $S_G$ of constant curvature $-1$. 
By construction, the natural cut locus structure of the image $x$ in $S_G$ of the center of $F_{2m}$ 
is precisely ${\cal C}$.
\hfill $\Box$

\bigskip

With a similar proof, one can show than every CL-structure on an arbitrary graph can be realized 
on a surface of constant curvature with at most $(n-p)$-singular points 
(i.e., on an Alexandrov surface with curvature bounded below, see  {\rm \cite{ST}} for the definition).
Here, $p$ is the number of vertices in $G$ of maximal degree.

\begin{ex}
{\rm The complete graphs $K_r$ and the multipartite graphs $K_{p_1,...,p_r}$
can be realized as cut loci on surfaces of constant curvature ($r, p_1,...,p_r \in I\!\!N$).

To obtain one realization of the Petersen graph as a cut locus, consider a regular $30$-gon $P$ in the hyperbolic plane $\mathbb{H}^2$ of constant curvature $-1$, with angles $2\pi/3$. Label the vertices of $P$, in circular order, by:
1, 2, 7, 9, 6, 1, 2, 3, 8, 10, 7, 2, 3, 4, 9, 6, 8, 3, 4, 5, 10, 7, 9, 4, 5, 1, 6, 8, 10, 5.
Now identify the edges having the same extremity labels, and get the desired surface $S$. Notice that $S$ is non-orientable.}
\end{ex}


\section{Stability}
\label{stability}

In this section we propose a notion of stability for cut locus structures, 
while in the next section we show that --roughly speaking-- 
stability is a generic property of CL-structures. For our goal, we need to further investigate the cyclic part of the cut locus; it was introduced and first studied by J. Itoh and T. Zamfirescu \cite{Itoh-Zamfirescu}.

The following result seems to be of some interest in its own right.

\begin{prop}
\label{cyclic_part}
The cyclic part of the cut locus depends continuously on the point; i.e.,

i) if $x_n \in S$, $x_n \to x$, and $y_n \in C^{cp}(x_n)$, $y_n \to y$, then $y \in C^{cp}(x)$, and

ii) if $x_n \in S$, $x_n \to x$, and $y \in C^{cp}(x)$, then there exist points $y_n \in C^{cp}(x_n)$ 
such that $y_n \to y$.
\end{prop}

\noindent{\sl Proof:} {\it i)}
It is well known that each limit of a sequence of geodesic segments is still a geodesic segment.
Assume we have two such sequences, say $\{\gamma_n\}_n$ and $\{\delta_n\}_n$, such that $\gamma_n$ and $\delta_n$ are both joining $x_n \in S$ to $y_n \in C^{cp}(x_n)$. Put $x_n \to x$, $y_n \to y$, and assume $\gamma_n \to \gamma$, $\delta_n \to \delta$.
Notice that $\gamma_n$ and $\delta_n$ determine a loop which is non homotopic to zero, because $y_n \in C^{cp}(x_n)$. So $\gamma \neq \delta$, and they also determine a loop non homotopic to zero; i.e., $y \in C^{cp}(x)$.

{\it ii)} For the second part, consider $x_n \in S$, $x_n \to x$.

The number $q$ of generating cycles in the cyclic part of a cut locus equals the first Betti number of $S$ \cite{sa}, hence it does not depend on the point in $S$. Therefore,
$$q(C^{cp}(x_n))=q(C^{cp}(x)).$$

Assume now that {\it (ii)} doesn't hold.
Then there exist a point $y \in C^{cp}(x)$ and a neighborhood $N_y \subset S$ such that $N_y \cap C^{cp}(x_n) = \emptyset$, for any $n$ sufficiently large.
Denote by $C^-$ the set of all such points $y$, and notice that $C^-$ is an open subset of $C^{cp}(x)$, with respect to the induced topology from $S$.

Notice that $C^{cp}(x_n)$ is a compact subset of $S$, hence $\lim_n C^{cp}(x_n)$ exists in the space of compact subsets of $S$, endowed with the usual topology induced by the Pompeiu-Hausdorff metric. 

Lemma \ref{qmn} and {\it (i)} show now that
$$q(C^{cp}(x))=\lim_n q(C^{cp}(x_n)) = q(\lim_n C^{cp}(x_n)) 
= q(C^{cp}(x) \setminus C^-) < q(C^{cp}(x)),$$
and a contradiction is obtained.
\hfill $\Box$

\begin{df}
{\rm Consider a CL-structure ${\cal C}$ on the graph $G$, a Riemannian surface $(S,g)$ and a point $x \in S$.
${\cal C}$ is called} {\it stable} {\rm with respect to $x$ in $S$ if}

i) {\rm $CLNS(x)={\cal C}$, and}

ii) {\rm there exists a neighborhood of $x$ in $S$, for all points $y$ of which $CLNS(y)={\cal C}$ holds.}
\end{df}

\begin{df}
{\rm The CL-structure ${\cal C}$ is called} {\it globally stable} {\rm if it is stable 
on all surfaces where it can be realized as a CLNS.}
\end{df}

Assume we have distinct pairs $(S, x)$ and $(S',x')$ of Riemannian surfaces $S$, $S'$ and points $x \in S$, $x' \in S'$ 
such that $CLNS(x)=CLNS(x')={\cal C}$. 
If ${\cal C}$ is stable with respect to $(S, x)$, it is not necessarily stable with respect to $(S',x')$,
as the following example shows.

\begin{figure*}
\label{Unstable}
\centering
  \includegraphics[width=1.0\textwidth]{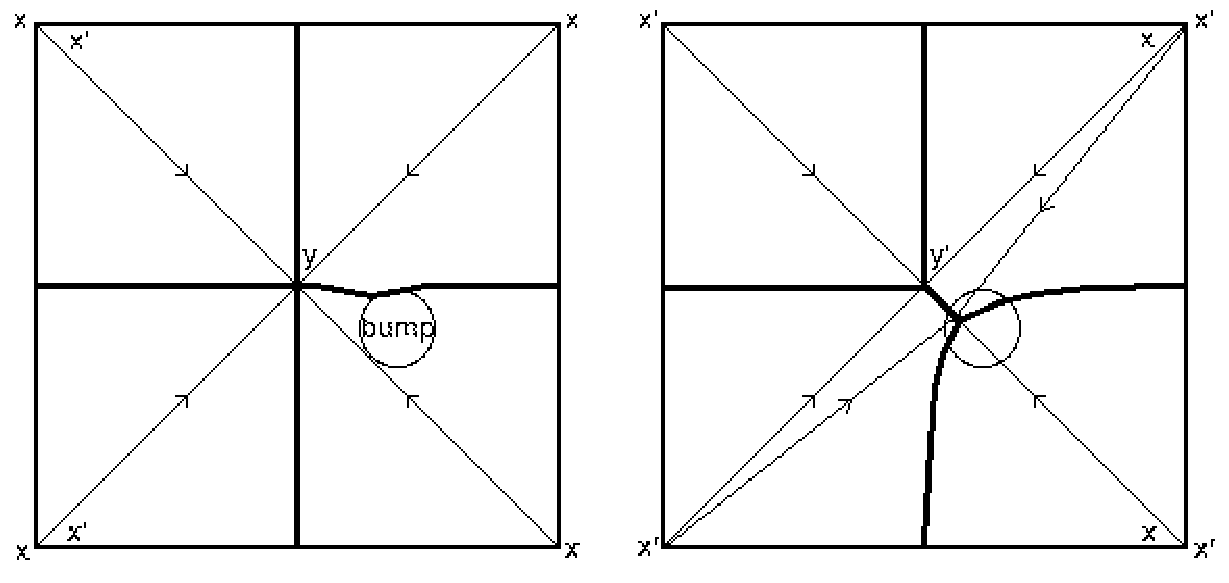}
\caption{Unstable cut locus structure.}
\end{figure*}

\begin{ex}
\label{bump-unstable}
i) {\rm Any CL-structure on a $k$-regular graph with $k >3$ is stable with res\-pect to 
the natural realization given by Theorem \ref{CL_ct_curv}.}

ii) {\rm We roughly explain here how to produce unstable CL-structures from those stable CL-structures at {\it (i)}.

Consider, for example, a square fundamental domain of a flat torus $T$ with a bump, see Figure 3 left.
The cut locus of the point $x \in T$, represented at the corners of the square, 
is the $4$-regular graph with one vertex $y$, as indicated by the thick line.
The four segments from $x$ to $y$ are also indicated with thin lines, and are not affected by the bump.
We choose $x$ such that one segment is tangent to the bump's boundary.

Now consider a point $x'$ arbitrarily close to $x$:
slightly move $x$ to ``the right'', for example, to $x'$, see Figure 3 right.
There we have only three segments from $x'$ to $y'$ (the center of figure with vertices at $x'$), 
those in the upper-left half-domain; they are all shorter than the geodesic joining $x'$ to $y'$ that crosses the bump, 
so $y'$ is a degree three vertex in $C(x')$.
There is another vertex of degree three in $C(x')$, also indicated in the figure 
together with the segments joining it to $x'$.
In this case, $C(x')$ is a $3$-regular graph with two vertices and two generating cycles.
J. Itoh and T. Sakai describe into details a similar procedure, see
Remark 2.7 in  {\rm \cite{IS}}.

In conclusion, the $4$-regular graph with one vertex is not stable with respect to $x$ in $T$.}
\end{ex}

\begin{thm}
\label{stability_3}
A cut locus structure on the graph $G$ is globally stable if and only if $G$ is a $3$-regular graph.
\end{thm}

\noindent{\sl Proof:}
Let ${\cal C}$ be a locus structure on $G$.

Assume first that $G$ is a $3$-regular graph; then its cyclic part is itself a $3$-regular graph.
Assume, moreover, that ${\cal C}$ is realized as ${\cal C}=CLNS(x)$, 
for some point $x$ on some Riemannian surface $S$.

Now, for points $x_n \in S$, $x_n \to x$, Proposition \ref{cyclic_part}
gives $\lim_n C^{cp}(x_n)=C^{cp}(x)$.

Assume that, for our sequence $\{x_n\}$, we have vertices $z_n$ in $C^{cp}(x_n)$ of degree $d$ larger than $3$, say $d=4$ (the case $d>4$ is similar).

Denote by $B_n^i$ the branches of $C^{cp}(x_n)$ incident to $z_n$; 
there exist segments $\gamma_n^i$, $\gamma_n^{'i}$ from $x_n$ to $z_n$, possibly with $\gamma_n^{i+1}=\gamma_n^{'i}$ ($i=1, ..., 4$, $\gamma_n^5=\gamma_n^1$) and a neighborhood $V_n$ of $z_n$ in $S$, such that one of the sectors around $z_n$ determined by $\gamma_n^i$, $\gamma_n^{'i}$ and $V_n$ contains $B_n^i \cap V_n$ but no other subsegment of a segment from $x_n$ to $z_n$.

Take some limit point $z$ of $z_n$; then $z\in C^{cp}(x)$, because $\lim_n C^{cp}(x_n)=C^{cp}(x)$, and $z$ has degree $3$ in $C^{cp}(x)$, by our assumption that $C^{cp}(x)$ is a cubic graph.
Therefore, there exists $1\leq i \leq 4$ such that the segments $\gamma_n^i$ and $\gamma_n^{'i}$ have a common limit $\gamma^i$, which is a segment from $x$ to $z$. Then, for $n$ large enough, $\gamma_n^i \cup \gamma_n^{'i}$ bounds a region of $S$ contractible to a point and intersecting $B_n^i \cap V_n \setminus \{z_n\}$.
Since $C^{cp}(x_n)$ intersects $\gamma_n^i \cup \gamma_n^{'i}$ only at $z_n$, it follows that $C^{cp}(x_n)$ contains a tree with the root at $z_n$, and a contradiction is obtained.

Concluding, the graph $C^{cp}(x_n)$ has to be cubic, and now 
$\lim_n C^{cp}(x_n)=C^{cp}(x)$ implies that the cyclic parts of $C(x)$ and $C(x_n)$ are isomorphic, 
and thus $G$ is stable.

\medskip

Assume now that $G$ is stable and it has a vertex $y$ of degree strictly larger than $3$, 
and consider a point $x$ in the Riemannian surface $S$ such that ${\cal C}=CLNS(x)$.
Then, by ``putting'' a bump tangent to one of the segments from $x$ to $y$
(i.e., modifying the metric on $S$ accordingly)
we obtain a new metric on $S$ with respect to which we still have ${\cal C}=CLNS(x)$,
but we have points $x'$ arbitrarily close to $x$ such that $CLNS(x') \neq {\cal C}$,
see Example \ref{bump-unstable} or Theorem \ref{global_stability}.
\hfill $\Box$

\bigskip

The following is, in some sense, opposite to Question \ref{all_x_CL}.

\begin{q}
\label{many_stable_CL}
{\rm How many stable CL-structures can exist on a given surface?}
\end{q}

Upper bounds on the number of cut locus structures on a graph are obtained in \cite{iv4}.


\section{Generic behavior}
\label{genericity}

We shall make use of the main result in \cite{Bu3}, given in the following as a lemma.
For, denote by ${\cal G}$ the {\sl space of all Riemannian metrics on the surface} $S$;
i.e., it is viewed as the space of sections of the bundle of positive definite symmetric matrices over $S$, 
endowed with the ${\cal C}^\infty$ Whitney topology \cite{Bu3}.

Recall that a metric $g$ on the surface $S$ is called {\sl cut locus stable} \cite{Bu3} if 
for any metric $h$ close to $g$ there is a diffeomorphism $\phi$ of the surface, depending continuously on $h$,
such that $\phi(C(x,g))=C(x,h)$; here, $C(x,g)$ denotes the cut locus of $x$ with respect to $g$.

\begin{lm}
\label{C(x)-stable} {\rm \cite{Bu3}}
For every point $x$ in $S$ there exists a set ${\cal B}_x$ of $C(x)$ stable metrics on $S$, open and dense in ${\cal G}$. 
Moreover, for any $g$ in ${\cal B}_x$, every ramification point of $C(x,g)$ is joined to $x$ by precisely three segments.
\end{lm}

In virtue of Definition \ref{diff_CL} and the remark following it, we can regard a CL-structure on the graph $G$ as a companion function $G \to \mathbb{Z}_2$.

A CL-structure is called {\it cubic} if its underlying graph is cubic.

\begin{thm}
\label{global_stability}
There exists an open and dense set in $S \times {\cal G}$, for every element $(x,g)$ of which the naturally defined cut locus structure $CLNS(x,g)$ is cubic and locally constant.
\end{thm}

\noindent{\sl Proof:}
Consider the subset ${\cal O}$ of $S \times {\cal G}$, containing all pairs $(x,g)$ for which the naturally defined cut locus structure $CLNS(x,g)$ is cubic.

The density of ${\cal O}$ in $S \times {\cal G}$ follows directly from Lemma \ref{C(x)-stable}.

We show next that ${\cal O}$ is open in $S \times {\cal G}$.
Assume this is not the case, hence there exist 
$(x,g) \in {\cal O}$ and a sequence $\{(x_n, g_n)\}_n \subset S \times {\cal G}$ convergent to $(x,g)$, such that $C^{cp}(x,g)$ is a cubic graph, but $C^{cp}(x_n,g_n)$ are not cubic.

For $n$ sufficiently large, $C^{cp}(x,g_n)$ are still cubic graphs, 
by Lemma \ref{C(x)-stable}. Moreover, an argument similar to the first part in the proof of Theorem \ref{stability_3} shows now that, 
for $g_n$ close enough to $g$, $C^{cp}(x,g_n)$ is a cubic graph homeomorphic to $C^{cp}(x,g)$.

Now, Theorem \ref{stability_3} shows that $C^{cp}(x_n,g_n)$ is a cubic graph homeomorphic to $C^{cp}(x,g_n)$, hence homeomorphic to $C^{cp}(x,g)$, and a contradiction is obtained.

Therefore, ${\cal O}$ is open in $S \times {\cal G}$ and, moreover, for every pair $(x,g)$ in ${\cal O}$ the naturally defined cut locus structure $CLNS(x,g)$ is locally constant.
\hfill $\Box$

\bigskip

The following result is well-known.

\begin{lm}
\label{graph_approx}
Every graph can be obtained from some cubic graph by edge contractions.
\end{lm}

Moving from a point with stable CL-structure to point with another stable CL-structure, one has to pass through a point with non-stable CL-structure, a CL-structure that --in particular-- lives on a non-cubic graph (see Theorems \ref{global_stability}, \ref{stability_3} and Lemma \ref{graph_approx}). At the level of CL-structures, one sees at a first step contraction(s) of one (or several) edge-strip(s), and at a second step ``blowing(s) up'' of all vertices of degree larger than $3$ to trees of order $3$. (A formal description is given in {\rm \cite{iv2}}.)

\bigskip

Non-isometric surfaces realizing the same graph $G$ as a cut locus (Theorem \ref{real_cl}) are homeomorphic to each other, since topologically they can be distinguished only by their genus, which is a function on the number of generating cycles of $G$.
Therefore, all distinct CL-structures on $G$ ``live'' on homeomorphic surfaces.
On the other hand, Theorem \ref{global_stability} shows in particular that equivalent CL-structures on $G$.


\bigskip

\noindent {\bf Acknowledgement } This work was begun during C. V\^\i lcu's stay at Kumamoto University, supported by JSPS; its completition at IMAR Bucharest was partially supported by the grant PN II Idei 1187.

The authors are indebted to all those persons who contributed with criticism, remarks and suggestions to the present form of the paper.


\bigskip

Jin-ichi Itoh

\noindent {\small Faculty of Education, Kumamoto
University
\\Kumamoto 860-8555, JAPAN
\\j-itoh@gpo.kumamoto-u.ac.jp}

\bigskip

Costin V\^\i lcu

\noindent {\small Institute of Mathematics ``Simion Stoilow'' of the
Romanian Academy
\\P.O. Box 1-764, Bucharest 014700, ROMANIA
\\Costin.Vilcu@imar.ro}

\end{document}